\documentclass{amsart}
\usepackage{color,graphicx,wrapfig}
\begin{document}

\title[Hidden dynamics and pulsating waves]{Hidden dynamics and the origin of
pulsating waves in Self-propagating High temperature Synthesis}
\author{R. Monneau}
\address{Ecole Nationale des Ponts et Chauss\'ees, CERMICS,
 6 et 8 avenue Blaise Pascal,
Cit\'e Descartes Champs-sur-Marne, 77455 Marne-la-Vall\'ee  Cedex 2, France}
\email{monneau@cermics.enpc.fr}
\urladdr{http://cermics.enpc.fr/\textasciitilde monneau/home.html}
\author{G.S. Weiss}
\address{Graduate School of Mathematical Sciences,
University of Tokyo,
3-8-1 Komaba, Meguro, Tokyo,
153-8914 Japan}
\email{gw@ms.u-tokyo.ac.jp}
\urladdr{http://www.ms.u-tokyo.ac.jp/\textasciitilde gw}
\thanks{G.S. Weiss has been partially supported by the Grant-in-Aid
15740100/18740086 of the Japanese Ministry of Education and partially supported
by a fellowship of the Max Planck Society. Both authors thank the Max Planck
Institute for Mathematics in the Sciences for the hospitality
during their stay in Leipzig.}
\keywords{SHS, solid combustion, singular limit, pulsating wave, 
high activation energy, Stefan problem, ill-posed, forward-backward, free boundary,
far from equilibrium dynamics}
\subjclass{80A25, 80A22, 35K55, 35R35}{}
\begin{abstract}
We derive the precise limit of SHS in the high
activation energy scaling suggested by B.J. Matkowksy-G.I. Sivashinsky 
in 1978 and by 
A. Bayliss-B.J. Matkowksy-A.P. Aldushin in 2002. 
In the time-increasing case the limit coincides with
the Stefan problem for supercooled water {\em with spatially
inhomogeneous coefficients}. In general it is a nonlinear forward-backward
parabolic equation {\em with discontinuous hysteresis term}.
\\
In the first part of our paper we give a complete characterization of the limit problem
in the case of one space dimension.
\\
In the second part we construct 
in any finite dimension
a rather large family
of pulsating waves for the limit problem.
\\
In the third part, we prove that for constant coefficients
the limit problem in any finite dimension
{\em does not admit non-trivial pulsating waves}.
\\
The combination of all three parts strongly suggests a relation
between the pulsating waves constructed in the present paper
and the numerically observed pulsating waves for finite activation energy
in dimension $n\ge 1$
and therefore provides a possible and surprising
explanation for the phenomena observed.
\\
All techniques in the present paper (with the exception
of the remark in the Appendix) belong to the
category far-from-equilibrium-analysis/far-from-bifurcation-point-analysis.
\end{abstract}

\maketitle

%%%%% private macros, f.e. the different environments
% the environment proof and proof* is defined automaticaly
\def\esssup{\hbox{\rm esssup }}
\def\pardist{\hbox{\rm pardist }}
\def\R{\hbox{\bf R}}
\def\Z{\hbox{\bf Z}}
\def\N{\hbox{\bf N}}
\def\sr{\hbox{\small\bf R}}
\def\sn{\hbox{\small\bf N}}
\def\supp{\hbox{\rm supp }}
\newtheorem{theo}{\textsc{Theorem}}[section]
\newtheorem{lem}[theo]{\textsc{Lemma}}
\newtheorem{pro}[theo]{\textsc{Proposition}}
\newtheorem{cor}[theo]{\textsc{Corollary}}
\newtheorem{defi}[theo]{\textsc{Definition}}
\newtheorem{rem}[theo]{\textsc{Remark}}
\newtheorem{hyp}[theo]{\textsc{Hypothesis}}
\newtheorem{pers}[theo]{\textsc{Perspectives}}
\newtheorem{conj}[theo]{\textsc{Conjecture}}

%%%% the main article
\section{Introduction}
The system
\begin{equation}
\begin{array}{l}\label{solid}
\partial_t u - \Delta u = v f(u)\\
\partial_t v = -v f(u)\; ,
\end{array}
\end{equation}
where $u$ is the normalized temperature, $v$ is the
normalized concentration of the reactant
and the non-negative nonlinearity $f$ describes the
reaction kinetics,
is a simple but widely used model for solid combustion
(i.e. the case of the Lewis number being $+\infty$).
In particular it is being used to model the industrial
process of Self-propagating High temperature Synthesis (SHS).
In the case of high activation energy interesting phenomena
like instability of planar waves, fingering and 
screw-like spin combustion waves are observed.\\
In \cite{1978} B.J. Matkowsky-G.I. Sivashinsky
derived from a special case of (\ref{solid}) a formal singular limit containing a jump condition
for the temperature on the interface.
Later 
it has been argued that the problem is for high
activation energy related to a Stefan problem describing
the freezing of supercooled water (see \cite{1978}, \cite[p. 57]{frankel1}).
Subsequently the Stefan problem for supercooled water became the basis
for numerous papers focusing on stability analysis of
(\ref{solid}), fingering, screw-like spin combustion waves etc.
(see for example \cite{frankel1},\cite{frankel2},\cite{frankel4},\cite{frankel5},\cite{frankel6},\cite{frankel7},\cite{hotspots}).
\\
Surprisingly there are few {\em mathematical} results on
the subject:
In \cite{loubeau} E. Logak-V. Loubeau proved existence of a planar
wave in one space dimension and gave a rigorous proof
for convergence as the activation energy goes to infinity.
\\
Instability of the planar wave
for a special linearization
(and high activation energy) is due to \cite{bonnet}.
\\[.5cm]
The present paper consists of three parts: in the first part we prove rigorously that 
in the case of one space dimension
the SHS system
converges to the irreversible
Stefan problem for supercooled water (cf. Theorem \ref{onespace}). In the time-increasing case 
we obtain also convergence in higher dimensions (see Theorem \ref{main}).
As the initial data of the reactant concentration
enter the equation as the activation energy goes to
infinity, our result also seems to provide a {\em possible explanation for the
numerically observed pulsating
waves} (cf. \cite{ikeda}, \cite{ikeda2} and \cite{hotspots}). 
\\
As a matter of fact, in the second part of  our paper (Theorem \ref{th:PTW}) we 
use the spatially inhomogeneous coefficients in order to
construct
a pulsating wave {\em for each periodic function} $v^0$ (or $Y^0$, respectively)
on $\R^n$, using the approach in 
\cite{hamel}. We also obtain the spin combustion waves (or ``helical waves'') on the cylinder mantle
(see Remark \ref{manifold}). 
\\
In contrast, 
we show in the third part (see Theorem \ref{th:LiouvillePTW}) 
that {\em for constant} $v^0$ in any finite dimension,
{\em no non-trivial pulsating waves exist}.
In addition, the formal result in the Appendix suggests that
in one space dimension the planar wave is stable.
\\
Taken together (cf. section \ref{conclusions}), 
our results strongly suggest a relation
between the pulsating waves constructed in the present paper
and the numerically observed pulsating waves for finite activation energy
in dimension $n\ge 1$ (cf. \cite{ikeda}, \cite{ikeda2} and \cite{hotspots})
and therefore provide a possible and surprising
explanation for the phenomena observed.
\\[.5cm]
In the original setting by B.J. Matkowsky-G.I. Sivashinsky \cite[equation (2)]{1978},
according to our result Theorem \ref{onespace},
\begin{equation}
\begin{array}{l}
\partial_t u_N - \Delta u_N = (1-\sigma_N)Ne^N v_N \exp(-N/u_N),\\
\partial_t v_N = -Ne^Nv_N \exp(-N/u_N),
\end{array}
\end{equation}
each limit $u_\infty$ of $u_N>0$ as $N\to\infty$
satisfies for $(\sigma_N)_{N\in \N} \subset  \subset [0,1)$
(for $\sigma_N \uparrow 1, N\uparrow\infty$ the limit in this scaling
is the solution of the heat equation; 
cf. Section \ref{appl1}
and Theorem \ref{onespace})
\begin{equation}\label{lim1}
\partial_t u_\infty-v^0 \partial_t \chi=\Delta u_\infty \hbox{ in } (0,+\infty)\times\Omega,
\end{equation}
where $v^0$ are the initial data of $v_\infty$ and
\[ \chi(t,x)\left\{\begin{array}{ll}
=0,& \esssup_{(0,t)} u_\infty(\cdot,x)<1\; ,\\
\in [0,1],& \esssup_{(0,t)} u_\infty(\cdot,x)=1\; ,\\
=1,& \esssup_{(0,t)} u_\infty(\cdot,x)>1\; .\end{array}\right.\]

In the SHS system with another scaling and a temperature threshold (see \cite[p. 109-110]{hotspots}),
\begin{equation}
\begin{array}{l}
\partial_t \theta_N - \Delta \theta_N\\
= (1-\sigma_N)
NY_N \exp((N(1-\sigma_N)(\theta_N-1))/
(\sigma_N + (1-\sigma_N)\theta_N))\chi_{\{\theta_N>\bar \theta\}},\\
\partial_t Y_N = -(1-\sigma_N)
NY_N \exp((N(1-\sigma_N)(\theta_N-1))/
(\sigma_N + (1-\sigma_N)\theta_N))\chi_{\{\theta_N>\bar \theta\}}
\end{array}
\end{equation}
where 
$N(1-\sigma_N)>>1, \sigma_N\in (0,1)$ and $\bar\theta \in (0,1)$,
each limit $\theta_\infty$ of $\theta_N$ satisfies (cf. Section \ref{appl2} and Theorem \ref{onespace})
\begin{equation}\label{lim2}
\partial_t \theta_\infty-Y^0 \partial_t \chi=\Delta \theta_\infty \hbox{ in } (0,+\infty)\times\Omega,
\end{equation}
where $Y^0$ are the initial data of $Y_\infty$ and
\[ \chi(t,x)\left\{\begin{array}{ll}
=0,& \esssup_{(0,t)} \theta_\infty(\cdot,x)<1\; ,\\
\in [0,1],& \esssup_{(0,t)} \theta_\infty(\cdot,x)=1\; ,\\
=1,& \esssup_{(0,t)} \theta_\infty(\cdot,x)>1\; .\end{array}\right.\]
To our knowledge this precise form of the limit problem, i.e. 
the equation with the discontinuous
hysteresis term,
has not been known. Even in the time-increasing case it does
{\em not} coincide with the formal result in \cite{1978}.
\\[.5cm]
In the case that $\theta_\infty$ (or $u_\infty$, respectively) is increasing
in time
and $v^0$ (or $Y^0$, respectively) is constant, 
our limit problem coincides with the
Stefan problem for supercooled water,
an extensively studied ill-posed problem (for a survey see \cite{dewynne}). 
As it is a forward-backward parabolic equation it is not
clear whether one should expect uniqueness (see \cite[Remark 7.2]{friedman} for an example
of non-uniqueness in a related problem).\\
On the positive side,
much more is known about the Stefan problem for supercooled water
than the SHS system, e.g.
existence of a finger (\cite{ivantsov}), instability of the finger (\cite{langer}),
one-phase solutions (\cite{friedman}) etc.;
those results, when combined with our convergence result,
suggest that similar properties should be true for the SHS
system.\\
It is interesting to observe that even in the 
time-increasing case our singular limit
{\em selects certain solutions} of the
Stefan problem for supercooled water.
For example, $u(t) =
(\kappa-1) \chi_{\{ t<1\}}+ \kappa \chi_{\{ t> 1\}}$
is for each $\kappa \in (0,1)$ a perfectly valid
solution of the Stefan problem for supercooled water,
but, as easily verified, it cannot be obtained from the ODE
\[ \partial_t u_\varepsilon(t) = -\partial_t\exp(-{1\over \varepsilon}\int_0^t \exp((1-1/(u_\varepsilon(s)+1))/\varepsilon)\> ds)\; \hbox{ as } \varepsilon\to 0\; .\]
Our approach does not involve stability or bifurcation analysis.
For showing the convergence we use standard compactness
and topological arguments. We prove that if the measure of
the ``burnt zone'' is small enough in a parabolic cylinder,
then in a cylinder of smaller radius there cannot be any
burnt part. For the construction of pulsating waves we use
the approach in \cite{hamel} as well as blow-up arguments,
Harnack inequality and so on. In order to obtain non-existence
of pulsating waves for constant $v^0$ we use a Liouville
technique.
\\
Let us conclude the introduction with a comparison to blow-up in semilinear
heat equations, as the main problem arising in our
convergence proof, i.e. excluding ``peaking of the solution'' or
burnt zones with very small measure, resembles
the blow-up phenomena in semilinear
heat equations. One could therefore hope to apply
methods used to exclude blow-up in low dimensions in order to
exclude peaking, say in two dimensions. There are
however problems:
First, here, we are dealing not with a single solution but
with the one-parameter family $u_\varepsilon$ concentrating at some ``peak''
as $\varepsilon$ gets smaller. Second, the
$\varepsilon$-problem is not a scalar equation but a 
degenerate system. Third, in contrast to blow-up, 
peaking would not necessarily imply $u_\varepsilon$
going to $+\infty$. 
Fourth, our limit problem is a two-phase problem
while most known results for blow-up
in semilinear
heat equations assume the solution to be non-negative.
Fifth, in our problem it does not make
much sense studying the onset of burning, say the 
first time when $u_\varepsilon\ge -\varepsilon$,
whereas studying the time of first blow-up can be very reasonable
for semilinear
heat equations.
The last and most important difference is that while
semilinear
heat equations are parabolic and therefore
well-posed in a sense, our limit problem contains a backward component
making it {\em ill-posed}.
\section{Notation}
Throughout this article $\R^n$ will be equipped with the Euclidean
inner product $x\cdot y$ and the induced norm $\vert x \vert\> .$
$B_r(x)$ will denote the open $n$-dimensional ball of center
$x\> ,$ radius $r$ and volume $r^n\> \omega_n\> .$
When the center is not specified, it is assumed to be $0.$\\
When considering a set $A\> ,$ $\chi_A$ shall stand for
the characteristic function of $A\> ,$
while
$\nu$ shall typically denote the outward
normal to a given boundary.
We will use the distance $\pardist$ with respect
to the parabolic metric
$d((t,x),(s,y))=\sqrt{|t-s|+|x-y|^2}$.\\
The operator $\partial_t$ will mean the partial derivative
of a function in the time direction, $\Delta$ the
Laplacian in the space variables and ${\mathcal L}^n$ the $n$-dimensional Lebesgue measure.\\
Finally ${\bf W}^{2,1}_p$ denotes the parabolic Sobolev space
as defined in \cite{nina}.
\section{Preliminaries}
In what follows, $\Omega$ is a bounded $C^1$-domain in $\R^n$ and
\[ u_\varepsilon\in \bigcap_{T \in (0,+\infty)}
{\bf W}^{2,1}_2((0,T)\times \Omega)\] is a strong solution
of the equation
\begin{equation}\label{sys}
\begin{array}{c}
\partial_t u_\varepsilon(t,x) - \Delta u_\varepsilon(t,x) = -v^0_\varepsilon(x)\partial_t\exp(-{1\over \varepsilon}\int_0^t g_\varepsilon(u_\varepsilon(s,x))\> ds)\; ,\\
u_\varepsilon(0,\cdot)=u_\varepsilon^0 \hbox{ in } \Omega, \nabla u_\varepsilon\cdot \nu=0 \hbox{ on }
(0,+\infty)\times \partial\Omega \; ;
\end{array}
\end{equation}
here $g_\varepsilon$ is a non-negative function
on $\R$ satisfying:\\
0) $g_\varepsilon$ is for each $\varepsilon\in (0,1)$
piecewise continuous with only one possible
jump at $z_0$, $g_\varepsilon(z_0-)=g_\varepsilon(z_0)=0$ in case of a jump, and
$g_\varepsilon$ satisfies for each $\varepsilon\in (0,1)$ and for every $z\in \R$ 
the bound $g_\varepsilon(z)\le C_\varepsilon (1+\vert z \vert)$.\\
1) $g_\varepsilon/\varepsilon \to 0$ as $\varepsilon\to 0$
on each compact subset of $(-\infty,0)$.\\
2) for each compact subset $K$ of $(0,+\infty)$
there is $c_K>0$ such that 
$\min(g_\varepsilon,c_K) \to c_K$ uniformly on $K$ as $\varepsilon\to 0$.\\
The initial data satisfy $0\le v^0_\varepsilon\le C<+\infty$, $v^0_\varepsilon$ converges in $L^1(\Omega)$ to
$v^0$ as $\varepsilon\to 0$,
$(u_\varepsilon^0)_{\varepsilon\in (0,1)}$ is bounded in $L^{2}(\Omega)$,
it is uniformly bounded from below by a constant $u_{\rm min}$,
and it converges in $L^1(\Omega)$ to $u^0$ as $\varepsilon\to 0$.
\begin{rem}
Assumption 0) guarantees existence of a global strong solution for each $\varepsilon\in (0,1)$.
\end{rem}
\section{The High Activation Energy Limit}\label{mainsec}
The following theorem has been proved in \cite{amuc}. Let us repeat the
statements and its proof for the sake of completeness.
\begin{theo}\label{main}
The family $(u_\varepsilon)_{\varepsilon\in (0,1)}$ is for each
$T\in (0,+\infty)$ precompact in $L^1((0,T)\times \Omega)$, and
each limit $u$ of $(u_\varepsilon)_{\varepsilon\in (0,1)}$ as a sequence
$\varepsilon_m\to 0,$ satisfies in the sense of distributions the 
initial-boundary value problem
\begin{equation}\label{stefan}
\partial_t u-v^0 \partial_t \chi =\Delta u \hbox{ in } (0,+\infty)\times\Omega,
\end{equation}
\[ u(0,\cdot)=u^0+v^0 H(u^0) \hbox{ in } \Omega\; , \; \nabla u\cdot \nu=0 \hbox{ on }
(0,+\infty)\times \partial\Omega\; ,\]
\[ \hbox{where } \chi(t,x) \left\{\begin{array}{ll}
\in [0,1],& \esssup_{(0,t)} u(\cdot,x)\le 0\; ,\\
=1,& \esssup_{(0,t)} u(\cdot,x)>0\; ,\end{array}\right.\]
and $H$ is the maximal monotone graph
\[ H(z) \left\{\begin{array}{ll}
=0,& z<0,\\ 
\in [0,1],& z=0,\\
=1,& z>0\; .\end{array}\right.\]
Moreover, $\chi$ is increasing in time and
$u$ is a supercaloric function.\\
If $(u_\varepsilon)_{\varepsilon\in (0,1)}$ satisfies
$\partial_t u_\varepsilon\ge 0$ in $(0,T)\times \Omega$,
then $u$ is a solution of the Stefan problem for
supercooled water, i.e.
\[ \partial_t u-v^0 \partial_t H(u) =\Delta u \hbox{ in } (0,+\infty)\times\Omega\; .\]
\end{theo}
\begin{rem}
Note that assumption 1) is only needed to prove the second statement
``If ....''. 
\end{rem}
\proof
{\bf Step 0 (Uniform Bound from below):}\\
Since $u_\varepsilon$ is supercaloric, it is bounded from
below by the constant $u_{\rm min}$.\\
{\bf Step 1 ($L^2((0,T)\times \Omega)$-Bound):}\\
The time-integrated function $v_\varepsilon(t,x) := \int_0^t u_\varepsilon(s,x)\> ds,$
satisfies
\begin{equation}
\partial_t v_\varepsilon(t,x) - \Delta v_\varepsilon(t,x) = w_\varepsilon(t,x) + u_\varepsilon^0(x)
\end{equation} where $w_\varepsilon$ is a measurable function satisfying
$0\le w_\varepsilon \le C.$
Consequently
\[ \int_0^T \int_\Omega (\partial_t v_\varepsilon)^2 \; + \;
{1\over 2} \int_\Omega |\nabla v_\varepsilon|^2(T)
\; = \; \int_0^T \int_\Omega (w_\varepsilon + u_\varepsilon^0) \partial_t v_\varepsilon
\]\[ \le \; {1\over 2} \int_0^T \int_\Omega (\partial_t v_\varepsilon)^2 \; + \;
{T\over 2} \int_\Omega (C+ |u_\varepsilon^0|)^2\; ,\]
implying
\begin{equation}\label{est1}
\int_0^T \int_\Omega u_\varepsilon^2 \; \le \;
T \int_\Omega (C+ |u_\varepsilon^0|)^2\; . 
\end{equation}
{\bf Step 2 ($L^2((0,T)\times \Omega)$-Bound for $\nabla\min(u_\varepsilon,M)$:}\\
For 
\[
G_M(z) := \left\{ \begin{array}{l}
z^2/2, z<M,\\
Mz-M^2/2, z\ge M\; ,\end{array}\right.\]
and
any $M\in\N,$
\[ \int_\Omega G_M(u_\varepsilon)-G_M(u^0_\varepsilon)
\; + \; \int_0^T \int_\Omega |\nabla \min(u_\varepsilon,M)|^2
\]\[ = \; \int_0^T \int_\Omega -v^0_\varepsilon \min(u_\varepsilon,M)
\partial_t \exp(-{1\over \varepsilon}\int_0^t g_\varepsilon(u_\varepsilon(s,x))\> ds)\; .\]
As $\partial_t \exp(-{1\over \varepsilon}\int_0^t g_\varepsilon(u_\varepsilon(s,x))\> ds)\le 0,$
we know that $\partial_t \exp(-{1\over \varepsilon}\int_0^t g_\varepsilon(u_\varepsilon(s,x))\> ds)$
is bounded in $L^\infty(\Omega;L^1((0,T))),$ and
\[ \int_0^T \int_\Omega -v^0_\varepsilon \min(u_\varepsilon,M)
\partial_t \exp(-{1\over \varepsilon}\int_0^t g_\varepsilon(u_\varepsilon(s,x))\> ds)\]\[
\le \; C\> \int_\Omega \sup_{(0,T)} \max(\min(u_\varepsilon,M),0) \;\le \;
CM {\mathcal L}^n(\Omega)\; .\]
\\
{\bf Step 3 (Compactness):}
\\
Let $\chi_M:\R\to\R$ be a smooth non-increasing function satisfying
$\chi_{(-\infty,M-1)}\le \chi_M \le \chi_{(-\infty,M)}$
and let $\Phi_M$ be the primitive such that $\Phi_M(z)=z$ for
$z\le M-1$ and $\Phi_M\le M$.
Moreover, let $(\phi_\delta)_{\delta\in (0,1)}$ be a family of mollifiers,
i.e.
$\phi_\delta\in C^{0,1}_0(\R^n;[0,+\infty))$ such that $\int \phi_\delta = 1$
and $\supp \phi_\delta \subset B_\delta(0)\> .$
Then, if we extend $u_\varepsilon$ and $v^0_\varepsilon$ by the value $0$ to the whole of $(0,+\infty)\times
\R^n,$ we obtain by the homogeneous Neumann data of $u_\varepsilon$ that
\[ \partial_t \left(\Phi_M(u_\varepsilon)*\phi_\delta\right)(t,x)
\]\[ = \; \left(\left(\chi_M(u_\varepsilon)\left(\chi_\Omega\Delta u_\varepsilon\> - \> v^0_\varepsilon
\partial_t \exp(-{1\over \varepsilon}\int_0^t g_\varepsilon(u_\varepsilon(s,x))\> ds)\right)
\right)*\phi_\delta\right)(t,x)
\]\[ = \; \int_{\sr^n}\chi_M(u_\varepsilon)(t,y)
\bigg( \chi_\Omega(y)\Delta u_\varepsilon(t,y)\]\[ - \> v^0_\varepsilon(y)
\partial_t \exp(-{1\over \varepsilon}\int_0^t g_\varepsilon(u_\varepsilon(s,y))\> ds)\bigg)
\phi_\delta(x-y)\> dy
\]\[ =\;\int_{\sr^n}\phi_\delta(x-y)\bigg(
-\chi_M'(u_\varepsilon(t,y))\chi_\Omega(y)|\nabla u_\varepsilon(t,y)|^2 \]\[- \>
\chi_M(u_\varepsilon(t,y)) v^0_\varepsilon(y)
\partial_t \exp(-{1\over \varepsilon}\int_0^t g_\varepsilon(u_\varepsilon(s,y))\> ds)\bigg)\]\[
+ \>\chi_M(u_\varepsilon(t,y)) \chi_\Omega(y)\nabla u_\varepsilon(t,y)
\cdot \nabla \phi_\delta(x-y)\> dy\; .\]
Consequently 
\[ \int_0^T \int_{\sr^n} |\partial_t \left(\Phi_M(u_\varepsilon)*\phi_\delta\right)|
\; \le \; C_1(\Omega,C,M,\delta,T)\]
and 
\[ \int_0^T \int_{\sr^n} |\nabla \left(\Phi_M(u_\varepsilon)*\phi_\delta\right)|
\; \le \; C_2(\Omega,M,\delta,T)\; .\]
It follows that $(\Phi_M(u_\varepsilon)*\phi_\delta)_{\varepsilon\in (0,1)}$
is for each $(M,\delta,T)$ precompact in $L^1((0,T)\times \R^n)$.
\\
On the other hand
\[ \int_0^T \int_{\sr^n} |\Phi_M(u_\varepsilon)*\phi_\delta-\Phi_M(u_\varepsilon)|
\; \le \; C_3 \left(\delta^2\int_0^T \int_\Omega |\nabla \Phi_M(u_\varepsilon)|^2\right)^{1\over 2}
\]\[ +\; 2(M-u_{\rm min})T{\mathcal L}^n(B_{\delta}(\partial\Omega))
\;\le \; C_4 (C,\Omega,u_{\rm min},M,T)\> \delta\; .\]
Combining this estimate
with the precompactness of $(\Phi_M(u_\varepsilon)*\phi_\delta)_{\varepsilon\in (0,1)}$ we obtain that $\Phi_M(u_\varepsilon)$ is for
each $(M,T)$ precompact in $L^1((0,T)\times \R^n)$.
Thus, by a diagonal sequence
argument, we may take a sequence $\varepsilon_m\to 0$ such that
$\Phi_M(u_{\varepsilon_m})\to z_M$ a.e. in $(0,+\infty)\times \R^n$
as $m\to\infty$, for every $M\in \N.$
At a.e. point of the set $\{ z_M < M-1\}, u_{\varepsilon_m}$
converges to $z_M.$ At each point $(t,x)$ of the remainder $\bigcap_{M\in \sn} 
\{ z_M \ge M-1\},$ the value $u_{\varepsilon_m}(t,x)$ must for large $m$
(depending on $(M,t,x)$) be larger than $M-2$. But that means that
on the set $\bigcap_{M\in \sn} 
\{ z_M \ge M-1\},$ the sequence $(u_{\varepsilon_m})_{m\in \sn}$ converges 
a.e. to $+\infty.$ It follows that $(u_{\varepsilon_m})_{m\in \sn}$
converges a.e. in $(0,+\infty)\times \Omega$ to a function 
$z:(0,+\infty)\times\Omega\to \R\cup\{ +\infty\}$.
But then, as $(u_{\varepsilon_m})_{m\in \sn}$ is for each $T\in (0,+\infty)$
bounded in $L^2((0,T)\times\Omega)$, $(u_{\varepsilon_m})_{m\in \sn}$
converges by Vitali's theorem (stating that a.e. convergence
and a non-concentration condition in $L^p$ imply 
in bounded domains $L^p$-convergence) for each $p\in [1,2)$ in $L^p((0,T)\times\Omega)$
to the weak $L^2$-limit $u$ of $(u_{\varepsilon_m})_{m\in \sn}$.
It follows that 
${\mathcal L}^{n+1}(\bigcap_{M\in \sn} 
\{ z_M \ge M-1\})={\mathcal L}^{n+1}(\{ u=+\infty\})=0$.
\\
{\bf Step 4 (Identification of the Limit Equation in $\esssup_{(0,t)} u>0$):}\\
Let us consider 
$(t,x)\in (0,+\infty)\times \Omega$ such that
$u_{\varepsilon_m}(s,x)\to u(s,x)$ for a.e. $s \in (0,t)$
and $u(\cdot,x) \in L^2((0,t))$.
In the case $\esssup_{(0,t)} u(\cdot,x)>0$,
we obtain by Egorov's theorem and assumption 2) that
$\exp(-{1\over {\varepsilon_m}}\int_0^t g_{\varepsilon_m}(u_{\varepsilon_m}(s,x))\> ds)\to 0$
as $m\to \infty$.
\\
{\bf Step 5 (The case $\partial_t u_\varepsilon \ge 0$):}\\
Let $(t,x)$ be such that $u_{\varepsilon_m}(t,x)\to u(t,x)=\lambda<0$: Then
by assumption 1),
\[ \exp(-{1\over {\varepsilon_m}}\int_0^t g_{\varepsilon_m}(u_{\varepsilon_m}(s,x))\> ds)
\ge \exp(-t{\max_{[u_{\rm min},\lambda/2]} g_{\varepsilon_m}\over {\varepsilon_m}})
\to 1\hbox{ as }m\to \infty\; .\]
\begin{rem}
We also obtain a rigorous convergence result
in the case of (higher dimensional) traveling waves 
with suitable conditions at infinity. In this case
our $L^2(W^{1,2})$-estimate (Step 2) implies
a no-concentration property of the time-derivative. 
\end{rem}
\section{Complete characterization of the limit equation in the case of one space dimension}
\begin{figure}
\begin{center}
\input{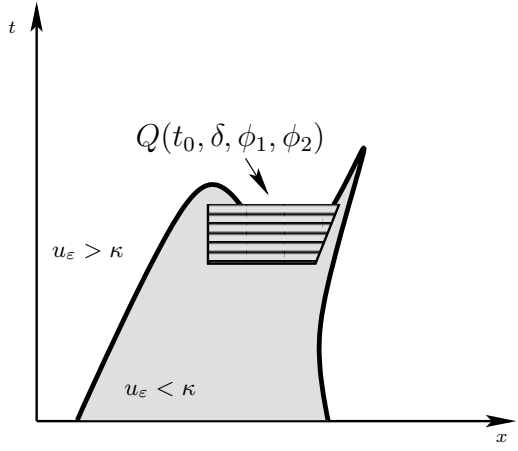}
\end{center}
\caption{Clearing out}\label{clearfig}
\end{figure}

\begin{figure}
\begin{center}
\input{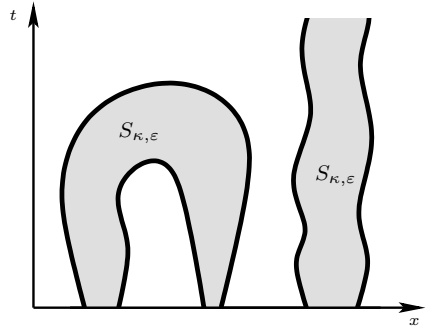}
\end{center}
\caption{The set $S_{\kappa,\varepsilon}$}\label{curvefig}
\end{figure}

\begin{figure}
\begin{center}
\input{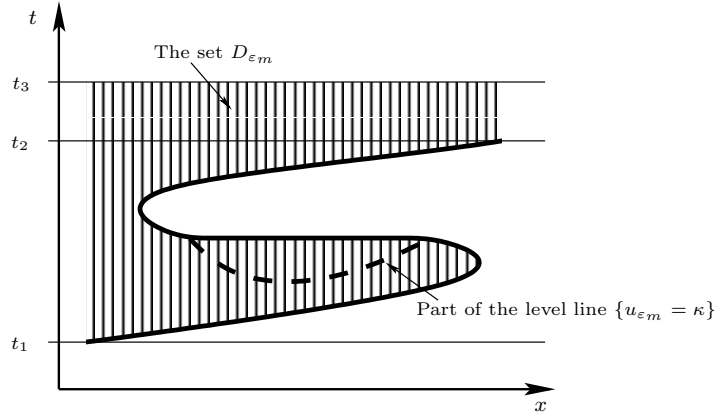}
\end{center}
\caption{The set $D_{\varepsilon_m}$}\label{D}
\end{figure}

The aim of this main section is the following theorem:
\begin{theo}\label{onespace}
Suppose in addition to the assumptions at the beginning of Section \ref{mainsec}
that the space dimension $n=1$ and that
the initial data $u_\varepsilon^0$ converge in $C^1$ to a function
$u^0$ satisfying $\nabla u^0\ne 0$ on $\{ u^0=0\}$.
Then the family $(u_\varepsilon)_{\varepsilon\in (0,1)}$ is for each
$T\in (0,+\infty)$ precompact in $L^1((0,T)\times \Omega)$,
and each limit $u$ of $(u_\varepsilon)_{\varepsilon\in (0,1)}$ as a sequence
$\varepsilon_m\to 0,$ satisfies in the sense of distributions the 
initial-boundary value problem
\begin{equation}
\partial_t u-v^0 \partial_t \chi =\Delta u \hbox{ in } (0,+\infty)\times\Omega,
\end{equation}
\[ u(0,\cdot)=u^0+v^0 H(u^0) \hbox{ in } \Omega\; , \; \nabla u\cdot \nu=0 \hbox{ on }
(0,+\infty)\times \partial\Omega\; ,\]
where $H$ is the maximal monotone graph
\[ H(z) \left\{\begin{array}{ll}
=0,& z<0,\\ 
\in [0,1],& z=0,\\
=1,& z>0\end{array}\right.\]
\[ \hbox{and } \chi(t,x) = H(\esssup_{(0,t)} u(\cdot,x)) \; \left\{\begin{array}{ll}
=0,& \esssup_{(0,t)} u(\cdot,x)< 0\; ,\\
\in [0,1],& \esssup_{(0,t)} u(\cdot,x)= 0\; ,\\
=1,& \esssup_{(0,t)} u(\cdot,x)>0\; ,\end{array}\right.\]
\end{theo}
Although we assume the space dimension from now on to be $1$,
we keep the multi-dimensional notation for the sake of
convenience. Moreover 
we extend $u_\varepsilon$ by even reflection at the lateral
boundary to a space-periodic solution on $[0,+\infty)\times \R$.\\
We start out with some elementary lemmata:
\\
\begin{lem}[Clearing out]\label{clearing}
There exists a continuous increasing function $\omega:[0,1)\to [0,+\infty)$
such that $\omega(0)=0$ and the following holds:
suppose that $\kappa<0$, that $\varepsilon \le \omega(\kappa)$, that $\delta \in (0,1)$ and that $u_\varepsilon \le (1+\omega(\delta))\kappa$ on the parabolic boundary
of the domain $Q(t_0,\delta,\phi_1,\phi_2) := \{ (t,x): 0\le t_0-2\delta<t<t_0, \phi_1(t)<x<\phi_2(t)\}$,
where $\phi_1<\phi_2$ are $C^1$-functions.
Then $u_\varepsilon \le \kappa$ in $Q(t_0,\delta,\phi_1,\phi_2)$ (cf. Figure \ref{clearfig}).
\end{lem}
\proof
Comparing $u_\varepsilon$ in $Q(t_0,\delta,\phi_1,\phi_2)$
to the solution of the ODE
\[ y'(t)=C g_\varepsilon(y)/\varepsilon, y(t_0-2\delta)=(1+\omega(\delta))\kappa\]
we obtain 
the statement of the lemma.
\begin{lem}\label{Scurves}
For almost all $\kappa<0$ the level set $((0,+\infty)\times\Omega)\cap\{u_\varepsilon=\kappa\}$
is a locally finite union of $C^1$-curves. For such $\kappa$
we define the set \[ 
S_{\kappa,\varepsilon} := \{ (t,x) \in (0,+\infty)\times\Omega): u_\varepsilon(t,x)>\kappa\]\[
\hbox{and there is no } (t_0,\delta,\phi_1,\phi_2) \in [0,+\infty)\times (0,1)
\times C^1\times C^1
\]\[\hbox{such that }
u_\varepsilon \le \kappa \hbox{ on the parabolic boundary
of the domain } Q(t_0,\delta,\phi_1,\phi_2)\}\]
(cf. Figure \ref{curvefig}).
Then
$\partial S_{\kappa,\varepsilon}=\bigcup_{j=1}^{N_{\kappa,\varepsilon}} \hbox{\rm graph } (g_{j,\kappa,\varepsilon})$
where $g_{j,\kappa,\varepsilon}: [0,T_{j,\kappa,\varepsilon}] \to \R$
are piecewise $C^1$-functions and $N_{\kappa,\varepsilon}$ is 
for small $\varepsilon$
bounded by a constant depending
only on the limit $u^0$ of the initial data.
\end{lem} 
\begin{rem}
For illustration of the definition of
$S_{\kappa,\varepsilon}$, imagine the set $\{ u_\varepsilon >\kappa\}$
filled with water in a $(t,x)$-plane where $t$ represents the height.
Our modification of $\{ u_\varepsilon >\kappa\}$ means then that
the water is now allowed to flow out through the ``bottom'' $\{ t=0\}$.
\end{rem}
{\sl Proof of Lemma \ref{Scurves}:}
By the definition of $S_{\kappa,\varepsilon}$ and by the fact that $u_\varepsilon$ is supercaloric,
each connected component of $\partial S_{\kappa,\varepsilon}$
is a piecewise $C^1$-curve and
touches $\{t=0\}$. Therefore the number of connected components
is for small $\varepsilon>0$ bounded by a constant
$\tilde N$ depending
only on the limit $u^0$ of the initial data.\\
Let us consider one connected component $\gamma$ of $\partial S_{\kappa,\varepsilon}$.
By  the definition of $S_{\kappa,\varepsilon}$ and by the fact that $u_\varepsilon$ is supercaloric, 
the derivative of the time-component of the piecewise $C^1$-curve $\gamma$
can change its sign at most once! Thus we can define for each curve $\gamma$
one or two piecewise $C^1$-functions of time such that
$\gamma$ is the union of the graphs of the two functions.
The total number of graphs $N_{\kappa,\varepsilon}$ is therefore bounded by $2\tilde N$.
\\[1cm]
{\sl Proof of Theorem \ref{onespace}:}
\newline
By Theorem \ref{main} we only have to prove that
$\chi=0$ in the set $\{\esssup_{(0,t)} u(\cdot,x)<0\}$. The main problem
is to exclude ``peaking'' of the solution $u_\varepsilon$, i.e.
tiny sets where $u_\varepsilon>\kappa$. Here we show that in the case of one
space dimension, ``peaking'' is not possible. More precisely, if
the measure of the set $u_\varepsilon>\kappa$ is small in a parabolic cube,
then $u_\varepsilon$ is strictly negative in the cube of half the radius, uniformly in $\varepsilon$.
The proof is carried out in two steps:
\newline
{\bf Step 1:}
Let $(\varepsilon_m)_{m\in \N}$ be the subsequence in the proof of Theorem \ref{main}.
As a.e. point $(t,x)\in ((0,+\infty)\times \R)\cap
\{ u<0\}$ is a Lebesgue point of the set $\{ u<0\}$,
we may assume that there exists $\kappa<0$
such that for any $\theta\in (0,1)$,
sufficiently small $r_0>0$ and every ${\varepsilon_m} \in (0,\varepsilon_0)$,
\[ {\mathcal L}^2(((t-2r_0,t)\times B_{2r_0}(x))\cap \{ u_{\varepsilon_m}<2\kappa\})
\ge \theta {\mathcal L}^2(((t-2r_0,t)\times B_{2r_0}(x)))\; .\]
{\bf Step 2:}
Suppose now that
$((t-r_0,t)\times B_{r_0}(x))\cap \{ u_{\varepsilon_m} > \kappa \}\ne\emptyset$
(where $\kappa$ is chosen such that $\{u_{\varepsilon_m}=\kappa\}$
and $\{u_{\varepsilon_m}=2\kappa\}$
are locally finite unions of $C^1$-curves):
then $((t-r_0,t)\times B_{r_0}(x))\cap \partial S_{\kappa,{\varepsilon_m}}$ and
$((t-r_0,t)\times B_{r_0}(x))\cap \partial S_{2\kappa,{\varepsilon_m}}$ must by Lemma \ref{Scurves} be 
connected to the parabolic boundary of
$(t-2r_0,t)\times B_{2r_0}(x)$.
The $L^2((0,T)\times \Omega)$-Bound 
\begin{figure}
\begin{center}
\input{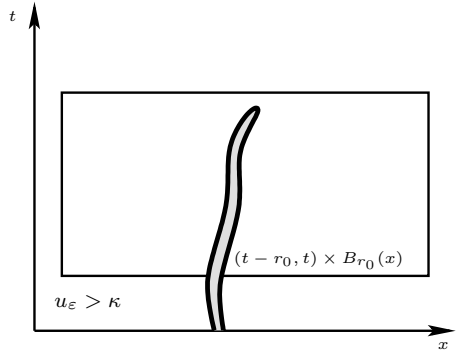}
\end{center}
\caption{Situation excluded by the $L^2(W^{1,2})$-estimate}\label{notposs1fig}
\end{figure}
for $\nabla\min(u_\varepsilon,M)$,
the fact that ${\mathcal L}^2(((t-2r_0,t)\times B_{2r_0}(x))\cap \{ u_{\varepsilon_m}<2\kappa\})
\ge \theta {\mathcal L}^2(((t-2r_0,t)\times B_{2r_0}(x)))$
and Lemma \ref{clearing}
imply now (see Figure \ref{notposs1fig}) that there must be an ``almost horizontal'' component
of $\partial S_{\kappa,{\varepsilon_m}}$ (cf. Figure \ref{horizontalfig}) with the following properties:
\begin{figure}
\begin{center}
\input{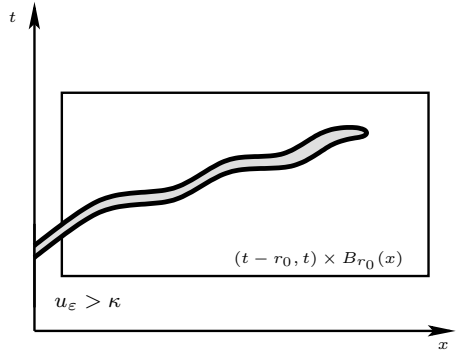}
\end{center}
\caption{The main task is to exclude almost horizontal propagation}\label{horizontalfig}
\end{figure}
\newline
for any $\delta\in (0,1)$, there are $t-r_0<t_1<t_2<t_3<t$ such that (see Figure \ref{D})
$t_3-t_1\to 0$ as ${\varepsilon_m} \to 0$, for some $j$
\[ |g_{j,\kappa,{\varepsilon_m}}(t_2)-g_{j,\kappa,{\varepsilon_m}}(t_1)|\ge c_1 >0\; ,\]
and
\[ {\mathcal L}^1(\{ y \in B_{r_0}(x) : u_{\varepsilon_m}(t_3,y) > 2\kappa\}) \le \delta\; ,\]
\[\int_{B_{r_0}(x)\cap\{ u_{\varepsilon_m}(t_3,\cdot) > 2\kappa\}}|u_{\varepsilon_m}(t_3,y)|\> dy
\; \le \delta\; .\]
We may assume that $c_1<r_0$, that
$g_{j,\kappa,{\varepsilon_m}}(t_2)=\sup_{(t_1,t_2)} g_{j,\kappa,{\varepsilon_m}}$,
that $g_{j,\kappa,{\varepsilon_m}}(t_2)>g_{j,\kappa,{\varepsilon_m}}(t_1)$
and that
$u_{\varepsilon_m}(s,y)>\kappa$ for some $d>0$ and $(s,y)\in (t_1,t_2)\times B_{r_0}(x)$ such that
$g_{j,\kappa,{\varepsilon_m}}(s)<y<d+g_{j,\kappa,{\varepsilon_m}}(s)$.
We define the set $D_{\varepsilon_m}:=$ \[ \{ (s,y): t_1<s<t_3, y < g_{j,\kappa,{\varepsilon_m}}(s)
\hbox{ for } s \in (t_1,t_2) \hbox{ and } y < g_{j,\kappa,{\varepsilon_m}}(t_2)
\hbox{ for } s \in [t_2,t_3)\}\]
(cf. Figure \ref{D}) and the cut-off function
$\phi(y) := \max(0,\min(y-g_{j,\kappa,{\varepsilon_m}}(t_1),g_{j,\kappa,{\varepsilon_m}}(t_2)-y))$.
It follows that
\[
c_1^2\kappa/4 + 2\delta + o(1)
\; \ge \;
o(1) + \int_{g_{j,\kappa,{\varepsilon_m}}(t_1)}^{g_{j,\kappa,{\varepsilon_m}}(t_2)} \phi(y) (u_{\varepsilon_m}(t_3,y)-\kappa)\> dy\]
\[ \ge \;
o(1) + \int_{g_{j,\kappa,{\varepsilon_m}}(t_1)}^{g_{j,\kappa,{\varepsilon_m}}(t_2)} \phi(y) u_{\varepsilon_m}(t_3,y)\> dy
- \kappa \int_{t_1}^{t_2} \phi(g_{j,\kappa,{\varepsilon_m}}(s)) g'_{j,\kappa,{\varepsilon_m}}(s)\> ds
\]\[ \ge \; 
\int_{g_{j,\kappa,{\varepsilon_m}}(t_1)}^{g_{j,\kappa,{\varepsilon_m}}(t_2)} \phi(y) u_{\varepsilon_m}(t_3,y)\> dy
- \int_{t_1}^{t_2} \phi(g_{j,\kappa,{\varepsilon_m}}(s)) u_{\varepsilon_m}(s,g_{j,\kappa,{\varepsilon_m}}(s)) g'_{j,\kappa,{\varepsilon_m}}(s)\> ds
\]\[ \ge \;
\int_{D_{\varepsilon_m}}
\phi \partial_t u_{\varepsilon_m}
\ge \; \int_{D_{\varepsilon_m}}
\phi \Delta u_{\varepsilon_m}
\]\[ = \; -\int_{D_{\varepsilon_m}}
\nabla \phi\cdot \nabla u_{\varepsilon_m}
\; + \; \int_{t_1}^{t_2} \phi(s,g_{j,\kappa,{\varepsilon_m}}(s)) \partial_x u_{\varepsilon_m}(s,g_{j,\kappa,{\varepsilon_m}}(s))
\chi_{\{ u_{\varepsilon_m}(s,g_{j,\kappa,{\varepsilon_m}}(s))=\kappa\}}
\> ds
\]\[ \ge \; -\int_{D_{\varepsilon_m}}
\nabla \phi\cdot \nabla u_{\varepsilon_m}
\; \to \; 0\hbox{ as }{\varepsilon_m}\to 0\; ,\]
a contradiction for small ${\varepsilon_m}$ provided that $\delta$ has been
chosen small enough; in the third inequality we used Lemma \ref{clearing},
and the convergence to $0$ is due to the uniform $L^2(W^{1,2})$-bound.
% \section{Strong Convergence}
% The result in the last section implies strong convergence of $\nabla u_{\varepsilon_m}$
% in $L^2$:
% \begin{pro}
% For $u_{\varepsilon_m},u$ as in Theorem \ref {onespace},
% \[ \nabla u_{\varepsilon_m} \to \nabla u\hbox{ strongly in } L^2((0,T);\Omega)\]
% for each $T\in (0,+\infty)$ as $m\to\infty$.
% \end{pro}
% {\sl Proof:}
% Multiplying the $\varepsilon_m$-equation by $\min(u_{\varepsilon_m}+\delta,0)$
% and integrating we obtain
% \[ \int_0^T \int_\Omega |\min(u_{\varepsilon_m}+\delta,0)|^2
% + {1\over 2}  \int_\Omega \min(u_{\varepsilon_m}(T)+\delta,0)^2 -  \min(u^0_{\varepsilon_m}+\delta,0)^2
% \le o(1)\]
% as $m\to\infty$. As 
% \[ \int_0^T \int_\Omega |\min(u+\delta,0)|^2
% + {1\over 2}  \int_\Omega \min(u(T)+\delta,0)^2 -  \min(u^0+\delta,0)^2 = 0\; ,\]
% we obtain
% \[ \limsup_{m\to\infty} \int_0^T \int_\Omega |\min(u_{\varepsilon_m}+\delta,0)|^2
% \le \int_0^T \int_\Omega |\min(u+\delta,0)|^2 \; .\]
\section{Applications}\label{applications}
Although the limit equation is an ill-posed problem,
the convergence to the limit seems to be robust
with respect to perturbations of the $\varepsilon$-system and the scaling:
here we mention two examples of different systems leading to
the same limit. Other examples can be found in mathematical
biology (see \cite{kawasaki} and \cite{maini}).\\
For the convergence results below we assume that
the space dimension is $1$.
\subsection{The Matkowsky-Sivashinsky scaling}\label{appl1}
We apply our result to the scaling in \cite[equation (2)]{1978}, i.e.
\begin{equation}
\begin{array}{l}
\partial_t u_N - \Delta u_N = (1-\sigma_N)Nv_N \exp(N(1-1/u_N)),\\
\partial_t v_N = -Nv_N \exp(N(1-1/u_N)),
\end{array}
\end{equation}
where the normalized temperature $u_N$ and the normalized
concentration $v_N$ are non-negative,
$(\sigma_N)_{N\in \N} \subset  \subset [0,1)$ (in the case
$\sigma_N \uparrow 1, N\uparrow\infty$ the limit equation
in the scaling as it is would be the heat equation,
but we could still apply our result to $u_N/(1-\sigma_N)$)
and the activation energy
$N\to \infty$.
\\
Setting $u_{\rm min} := -1, \varepsilon := 1/N, u_\varepsilon := u_N-1$ and
\[ g_\varepsilon(z) := \left\{\begin{array}{ll}
\exp((1-1/(z+1))/\varepsilon), z>-1\\
0, z\le -1\end{array}\right .\]
and integrating the equation for $v_N$ in time, 
we see that the assumptions of Theorem \ref{onespace} are satisfied
and we obtain
that each limit $u_\infty,\sigma_\infty$ of $u_N,\sigma_N$ satisfies
\begin{equation}\label{matsiv}
\partial_t u_\infty-(1-\sigma_\infty)
v^0 \partial_t H(\esssup_{(0,t)} u_\infty)=\Delta u_\infty \hbox{ in } (0,+\infty)\times\Omega,
\end{equation}
\[ u_\infty(0,\cdot)=u^0+v^0H(u^0) \hbox{ in } \Omega\; , \; \nabla u_\infty\cdot \nu=0 \hbox{ on }
(0,+\infty)\times \partial\Omega\; ,\]
where $v^0$ are the initial data of $v_\infty$.
Moreover, $\chi$ is increasing in time and
$u_\infty$ is a supercaloric function.
\subsection{SHS in another scaling with temperature threshold}\label{appl2}
Here we consider (cf. \cite[p. 109-110]{hotspots}), 
i.e.
\begin{equation}
\begin{array}{l}
\partial_t \theta_N - \Delta \theta_N \\
= (1-\sigma_N)
NY_N \exp((N(1-\sigma_N)(\theta_N-1))/
(\sigma_N + (1-\sigma_N)\theta_N))\chi_{\{\theta_N>\bar \theta\}},\\
\partial_t Y_N = -(1-\sigma_N)
NY_N \exp((N(1-\sigma_N)(\theta_N-1))/
(\sigma_N + (1-\sigma_N)\theta_N))\chi_{\{\theta_N>\bar \theta\}}
\end{array}
\end{equation}
where 
$N(1-\sigma_N)>>1, \sigma_N\in (0,1)$ and
the constant $\bar\theta \in (0,1)$ is a threshold parameter
at which the reaction sets in.
\\
Setting $u_{\rm  min}=-1, \varepsilon := 1/(N(1-\sigma_N)), 
\kappa(\varepsilon):= 1-\sigma_N, u_\varepsilon := \theta_N-1$,
\[ g_\varepsilon(z) := \left\{\begin{array}{ll}
\exp((z/(\kappa(\varepsilon)z+1))/\varepsilon), z>\bar\theta-1\\
0, z\le \bar\theta-1\end{array}\right .\]
and integrating the equation for $Y_N$ in time,
we see that the assumptions of Theorem \ref{onespace} are satisfied
and we obtain that each limit $u_\infty$ of $u_N$ satisfies
\begin{equation}\label{threshold}
\partial_t u_\infty-v^0 \partial_t H(\esssup_{(0,t)} u_\infty)=\Delta u_\infty \hbox{ in } (0,+\infty)\times\Omega,
\end{equation}
\[ u_\infty(0,\cdot)=u^0+v^0H(u^0) \hbox{ in } \Omega\; , \; \nabla u_\infty\cdot \nu=0 \hbox{ on }
(0,+\infty)\times \partial\Omega\; ,\]
where $v^0$ are the initial data of $v_\infty$.
Moreover, $\chi$ is increasing in time and
$u_\infty$ is a supercaloric function.
\section{Existence of Pulsating Waves}\label{expw}
The aim of this section is to construct pulsating waves
for the limit problem. 
For the sake of clarity we have chosen not to present the most general
result in the following theorem. Moreover we confine ourselves
to the one-phase case.

\begin{theo}\label{th:PTW}{\bf (Existence of pulsating waves)}\\
Let us consider a H\"older continuous function $v^0$ defined on $\R^n$ that
satisfies
$$
%\begin{equation}
\label{eq:v0}
v^0(x)\ge 1 \quad \mbox{and}\quad v^0(x+k)=v^0(x) \quad\mbox{for every}\quad
k\in \Z^n,\quad x\in\R^n\; .   
%\end{equation}
$$
Given a unit vector $e\in \R^n$ and a velocity $c>0$, 
there exists a solution $u(t,x)$
of the one-phase problem
\begin{equation}\label{eq:uv0}
\left\{\begin{array}{l}
\partial_t u -v^0\partial_t \chi_{\left\{u\ge 0\right\}}=\Delta u
\quad\mbox{on}\quad \R\times \R^n\\
\\
\partial_{t}u\ge 0\quad\mbox{and}\quad -\mu_0:= -\int_{[0,1)^n} v^0\le u\le 0\; ,
\end{array}\right.
\end{equation}
which satisfies
\begin{equation}\label{eq:uv01}
\left\{\begin{array}{lcl}
\displaystyle{u(t,x+k)=u(t-\frac{e\cdot k}{c},x)
\quad\mbox{for every}\quad
k\in \Z^n,\quad (t,x)\in\R\times \R^n}\\
\\
u(t,x)=0 \quad\mbox{for}\quad x\cdot e -ct \le 0\quad\mbox{and}\quad
\displaystyle{\limsup_{x\cdot e -ct \to +\infty } u(t,x) = -\mu_0}\; ,
\end{array}\right.
\end{equation}
where the last limit is uniform as  $x\cdot e -ct$ tends to $+\infty$.
\end{theo}
\begin{rem}\label{manifold}
By modifications of the following proofs and of the theory in
\cite{hamel}, it is possible to replace $\R^n$ in Theorem \ref{th:PTW}
by a smooth source manifold. Taking for example
$S^1\times \R$, we obtain the screw-like pulsating waves observed
in spin combustion (also called ``helical waves''; see for example \cite{ikeda}, \cite{ikeda2}, \cite{hotspots},\cite{ikeda}). 
\end{rem}
Let us transform the problem by the so-called Duvaut
transform (see \cite{rodrig}), setting 
$w(t,x)=-\int_t^{+\infty} u(s,x)\ ds$. 
In this section we will prove the existence of
a pulsating wave $w$.
More precisely, Theorem \ref{th:PTW} is a corollary of
of the following result which will be proved later.
\begin{theo}\label{pro:PTWOP}{\bf (Pulsating  waves for the obstacle problem)}\\
Under the assumptions of Theorem \ref{th:PTW}, there exists a function 
$w(t,x)$ solving the obstacle problem
\begin{equation}\label{eq:w}
\left\{\begin{array}{l}
\partial_t w =\Delta w-v^0\ \chi_{\left\{w> 0\right\}}
\quad\mbox{on}\quad \R\times \R^n\; ,\\
\\
w\ge 0 \quad , \quad -\mu_0\le \partial_t w\le 0, 
\quad \partial_{tt}w\ge 0\; ,
\end{array}\right.
\end{equation}
with the conditions
\begin{equation}\label{eq:w1}
\left\{\begin{array}{l}
\displaystyle{w(t,x+k)=w(t-\frac{e\cdot k}{c},x)
\quad\mbox{for every}\quad
k\in \Z^n,\quad (t,x)\in\R\times \R^n}\; ,\\
\\
w(t,x)=0 \;\mbox{ for }\; x\cdot e -ct \le 0
\quad\mbox{and}\quad
\partial_t w(t,x) \longrightarrow - \mu_0 \;\mbox{ as }\; x\cdot e -ct
\longrightarrow +\infty\; .
\end{array}\right.
\end{equation}
The convergence is uniform as $x\cdot e -ct$ tends to $+\infty$.
\end{theo}

%%%%%%%%%%%%%%%%%%%%%%%%%%%%%%%%%%%%%%%%%%%%%%%%%%%%%%%%%%%%%%%%%%%%%%%%%%%%%%%%
\noindent {\sl Proof of Theorem \ref{th:PTW}}\\
Simply set $u(t,x):=\partial_tw(t,x)$ with $w$ given by Theorem
\ref{pro:PTWOP}, and use the fact that
%the Lebesgue measure of the free boundary $\partial\left\{w>0\right\}$ is
%zero, and 
$\chi_{\left\{u<0\right\}}=\chi_{\left\{w>0\right\}}$. To check this last
property, 
it is sufficient to exclude the case where $w(t_0,x_0)>0$ and $\partial_t
w(t_0,x_0)=0$ at some point $(t_0,x_0)$: using the fact that
$\partial_t w$ is caloric in
$\left\{w>0\right\}$ as well as the strong maximum principle, we deduce that
$\partial_t w(t,x)=0$ for $t\in(-\infty, t_0]$ and $x$ in a neighborhood of
$x_0$. This contradicts the last line of (\ref{eq:w1}).\\

%%%%%%%%%%%%%%%%%%%%%%%%%%%%%%%%%%%%%%%%%%%%%%%%%%%%%%%%%%%%%%%%%%%%%%%%%%%%%%%%
\noindent {\sl Proof of Theorem \ref{pro:PTWOP}}\\
We will prove the existence of an unbounded solution $w$ in six steps,
approximating $w$ by bounded solutions of a truncated equation, for which
we can apply the existence of pulsating fronts due to Berestycki and
Hamel \cite{hamel}.\\
%%%%%%%%%%%%%%%%%%%%%%%%%%%%%%%%%%%%%%%%%%%%%%%%%%%%%%%%%%%%%%%%%%%%%%%%%%%%%%%%

\noindent {\bf Step 1: Approximation by bounded solutions and estimates of
  the velocity}\\
For any $0<A<M$, let us start by approximating the function
$\chi_{(0,+\infty )}$ by the characteristic function $g=\chi_{(0,A]}$. 
In that case we can compute explicitly the {\em traveling wave} $(\phi,c_0)$ (unique up to
translations of $\phi$) of
\begin{equation}\label{eq:1Dpb}
c_0\phi'=\phi''-g(\phi), \quad  \phi'\le 0\quad
\mbox{on}\quad \R,\quad \phi(-\infty )=M\quad \mbox{and}\quad  \phi(+\infty )=0\quad :  
\end{equation}
Let us define for $c_0>0, M>0,s_0\in (-\infty,0)$ and 
$s_1\in (s_0,0)$
$$\phi(s)=\left\{\begin{array}{ll}
M\left(1-e^{c_0(s-s_0)}\right) &\quad \mbox{for}\quad s\in (-\infty ,s_1]\; ,\\
\displaystyle{\frac{1}{c_0^2}\left(e^{c_0s}-1-c_0s\right)} &\quad
\mbox{for}\quad  s\in [s_1,0]\; ,\\
 0 &\quad \mbox{for}\quad  s\in [0,+\infty )\; .
\end{array}\right.$$
For any $A\in (0,M)$ and for suitable $s_0,s_1$, we see that $\phi$ is continuous and satisfies $\phi(s_1)=A$,
which fixes the parameter $s_1$ as a function of $A$. Moreover
we see that $\phi$ is of class $C^1$ if and only if $s_1=-c_0M$ and
\begin{equation}\label{eq:c_0}
M-A=\frac{1}{c_0^2}\left(1-e^{-c_0^2M}\right)\; .
\end{equation}
Thus $A$ is determined in terms of the velocity $c_0$ and $M$.
The above calculations show in particular that $\phi(c_0t-e\cdot x)$ is a good bounded
approximation of the solution of (\ref{eq:w}) -(\ref{eq:w1}) in the case $v^0=1$,
i.e. the traveling wave case.\\
In all that follows let $A$ be given by (\ref{eq:c_0}).
\\
Now, when the function $g$ in (\ref{eq:1Dpb}) is replaced by a 
Lipschitz continuous function whose support is a compact interval,
there are known results on the existence of pulsating waves.
For such $g$, it is possible to apply Theorem 1.13 of
Berestycki, Hamel \cite{hamel}, which states the existence (and uniqueness up
to translation in time) of bounded pulsating
solutions traveling at a unique velocity.
Bearing that in mind, we define $g_M$ as a Lipschitz regularization of the characteristic
function
$g$ such that $\mbox{supp}\ g_M= [0,A]$ and -- for later use --
\begin{equation}\label{eq:g'g}
\begin{array}{l}
g_M=1 \quad\mbox{on} \quad [1/M,A/2],\quad 0\le g_M\le 1\quad\mbox{on}\quad \R\; ,\\
\mbox{and}\quad \quad g_M'\ge 0, \quad g_M''\le 0 \quad\mbox{on}\quad (0,A/2)\; .
\end{array}
\end{equation}
Let us call $c^0_M$ the unique velocity of the traveling wave equation
(\ref{eq:1Dpb}) with $g$ replaced by $g_M$.
As $(\phi,c_0)$ can be shown to be unique up to
translations of $\phi$), $c^0_M\to c_0$ as $g_M\to g$.
\\
Then there exists by Theorem 1.13 of
\cite{hamel} a bounded pulsating wave 
$w_M$ traveling at velocity $c_M$ such that
$$\left\{\begin{array}{l}
\partial_t w_M =\Delta w_M -v^0g_M(w_M) 
\quad\mbox{on}\quad \R\times \R^n\; ,\\
\\
\partial_t w_M\le 0\quad\mbox{and}\quad \displaystyle{\limsup_{x\cdot e -c_M t \to -\infty } w_M(t,x) = 0}\le w_M
\le M= \displaystyle{\liminf_{x\cdot e -c_M t \to +\infty } w_M(t,x)}\\
\\
\displaystyle{\mbox{and}\quad w_M(t,x+k)=w_M(t-\frac{e\cdot k}{c_M},x)
\quad\mbox{for every}\quad
k\in \Z^n,\quad (t,x)\in\R\times \R^n\; .}
\end{array}\right.$$
From the assumption $1\le v^0$ and the comparison principle (see Lemma 3.2 and
3.4 of \cite{hamel}) we infer that
the velocities 
satisfy the ordering
$c^0_M\le c_M$. Similarly, defining $\lambda =
||v^0||_{L^\infty}$ and comparing to $w_M(\lambda t,\sqrt{\lambda}x)$, we
get $c_M\le c^0_M \sqrt{\lambda}$.\\
Furthermore the above comparison principles tell us that the velocity $c_M$
(resp. $c^0_M$) is
continuous and non-decreasing in $M$.\\

For all that follows let $c>0$ be an arbitrary but fixed velocity
for which we want to construct the pulsating wave.
Then, 
for any $M>0$ we can adjust $A\in (0,M)$ such that
for $c_0$ defined by (\ref{eq:c_0}),
$$c=c_M \in [c_0/2,2c_0\sqrt{||v^0||_{L^\infty}}]\; .$$

In order to pass to the limit as $M\to +\infty$, we need to get some
bounds on the solution first. To this end, rotating space-time proves to be
very convenient:

%%%%%%%%%%%%%%%%%%%%%%%%%%%%%%%%%%%%%%%%%%%%%%%%%%%%%%%%%%%%%%%%%%%%%%%%%%%%%%%%
\noindent {\bf Step 2: Space-time transformation and first estimates on the time derivatives}\\
Let us introduce the function $\tilde{w}_M$ defined by 
 $\tilde{w}_M(s,x)= w_M(\frac{s+e\cdot x}{c},x)$ which is periodic in $x$ and
 satisfies
 $$L\tilde{w}_M = v^0(x)g_M(\tilde{w}_M)\quad\mbox{with}\quad L\tilde{w}_M =\Delta \tilde{w}_M + \partial_{ss}\tilde{w}_M - 2
 \partial_{e,s}\tilde{w}_M-c\partial_s\tilde{w}_M$$
 and 
$\lim_{s\to -\infty }\tilde{w}_M(s,x) = M$, $\lim_{s\to
   +\infty }\tilde{w}_M(s,x) = 0$ uniformly with respect to $x$.
Using (\ref{eq:g'g}), we obtain $L\partial_{s}\tilde{w}_M-c_1\partial_{s}\tilde{w}_M= 0$ and
 $L\partial_{ss}\tilde{w}_M-c_1\partial_{ss}\tilde{w}_M\le 0$ on
 $\left\{\tilde{w}_M < A/2\right\}$ for
 $c_1(s,x)=v^0(x)g_M'(\tilde{w}_M(s,x))\ge 0$ on this set.
We deduce from the maximum principle (see Lemma 3.2 and
3.4 of \cite{hamel}) that
for any $s_0\in \R$ such that $\sup_{[s_0,+\infty )\times \R^n}\tilde{w}_M\le
 A/2$,
\begin{equation}\label{eq:td}
\begin{array}{l}
\inf_{[s_0,+\infty )\times \R^n} \partial_s\tilde{w}_M = \inf_{\R^n}
\partial_s\tilde{w}_M (s_0,\cdot)\\
\mbox{and}\quad 
\inf_{[s_0,+\infty)\times \R^n} \partial_{ss}\tilde{w}_M = \min\left(0,\ 
\inf_{\R^n}\partial_{ss}\tilde{w}_M (s_0,\cdot)\right) \; .
\end{array}
\end{equation}
%%%%%%%%%%%%%%%%%%%%%%%%%%%%%%%%%%%%%%%%%%%%%%%%%%%%%%%%%%%%%%%%%%%%%%%%%%%%%%%%
\noindent {\bf Step 3: Bound of the solution from above}\\
  From the fact that $g_M$ is bounded by $1$, and from the Harnack inequality, we deduce that there exists a constant $C_H\in (1,+\infty)$ such
  that for any $r>0$ and for any point $(t_0,x_0)$
\begin{equation}\label{eq:Harnack}
\sup_{B_r(x_0)}w_M(t_0-r^2,\cdot )  \le C_H\left(\inf_{B_r(x_0)}w_M(t_0,\cdot)
    + r^2 \lambda \right) 
\end{equation}
where $\lambda=||v^0||_{L^\infty}$. 
For $\tilde{w}_M$ that means  that -- setting $s_0=ct_0-e\cdot x_0$ --
   $$\sup_{y\in B_{\sqrt{n}/2}(0)} \tilde{w}_M(s_0-cr^2 -e\cdot y,x_0+y) \le
   C_H\left(\inf_{y\in B_{\sqrt{n}/2}(0)} \tilde{w}_M(s_0-e\cdot y,x_0+y) 
     + r^2 \lambda \right)$$
for $r\ge \sqrt{n}/2$.
We will now use the fact that the unit cell $(-1/2,1/2)^n$ is contained
in the ball $B_{\sqrt{n}/2}(0)$. 
Using first the monotonicity of $\tilde{w}_M$ in the variable $s$, and
second the periodicity of $\tilde{w}_M(\tau,y)$ in $y$, we get for 
$\tau_0:=s_0-\sqrt{n}/2$
\begin{equation}\label{eq:ns1}
\displaystyle{\sup_{\R^n} \tilde{w}_M(\tau_0-cr^2+\sqrt{n},\cdot)}
    \le \displaystyle{C_H\left(\inf_{\R^n} \tilde{w}_M(\tau_0,\cdot)
     + r^2 \lambda \right)}\; .
\end{equation}
By a translation in time we may assume that
\begin{equation}\label{eq:ori}
 0=\inf\left\{\tau:\quad \tilde{w}_M(s,x) \le  1/M
   \quad\mbox{for}\quad s\ge\tau,\quad x\in\R^n\right\}  
\end{equation}
and get the bound
\begin{equation}\label{eq:bound}
\tilde{w}_M(s,x)\le \max\left(1/M, \alpha -\beta s\right)
\end{equation}
for some constants $\alpha,\beta\in (0,+\infty)$ and every large positive $M$.\\

%%%%%%%%%%%%%%%%%%%%%%%%%%%%%%%%%%%%%%%%%%%%%%%%%%%%%%%%%%%%%%%%%%%%%%%%%%%%%%%%
\noindent {\bf Step 4: Passing to the limit}\\
By estimate (\ref{eq:bound}), we can pass to the limit 
as $M\to +\infty$ and obtain
$M-A\to 1/c_0^2$. Moreover,
passing to a subsequence if necessary,
$\tilde{w}_M$ converges locally in ${\bf W}^{2,1}_p$ to $\tilde{w}$ satisfying
$$\left\{\begin{array}{l}
\partial_t w\le 0\quad\mbox{and}\quad \displaystyle{\limsup_{x\cdot e -c t \to -\infty } w(t,x) = 0}\le w
\\
\displaystyle{\mbox{and}\quad w(t,x+k)=w(t-\frac{e\cdot k}{c},x)
\quad\mbox{for every}\quad
k\in \Z^n,\quad (t,x)\in\R\times \R^n\; .}
\end{array}\right.$$

Furthermore,
we obtain for $w$ related to $\tilde{w}$ by
$\tilde{w}(s,x)=w\left(\frac{s+e\cdot x}{c},x\right)$
that
$$w_t=\Delta w -v^0 \chi_{\left\{w>0\right\}}\; ;$$
here we used the fact that
$w$, being locally a ${\bf W}^{2,1}_p$-function, satisfies
$\partial_t w=0=\Delta w$ a.e. on the set $\{ w=0\}$.

In order to conclude $\partial_{tt} w\ge 0$ the following
non-degeneracy property will prove to be necessary:

%%%%%%%%%%%%%%%%%%%%%%%%%%%%%%%%%%%%%%%%%%%%%%%%%%%%%%%%%%%%%%%%%%%%%%%%%%%%%%%%
\noindent {\bf Step 5: Non-degeneracy property and bound from below}\\
Let us assume that $w_M(t_0,x_0)\in (1/M,A/2)$. Using the fact
that $v^0(x)g_M(z)\ge 1$ for $z\in [1/M,A/2]$, we can use the
usual parabolic maximum principle, comparing $\max(w_M,1/M)$ to
the function $$h(t,x)=w_M(t_0,x_0) + \frac{1}{4n}|x-x_0|^2 +\frac{1}{4n}(t_0-t)$$ on the set
$$\left\{1/M <  w_M< A/2\right\}\cap Q^-_r(t_0,x_0)\; ,$$ where
$Q^-_r(t_0,x_0)=\left\{(t,x): \quad t_0-r^2\le t\le t_0,\quad |x-x_0|\le
  r\right\}$. 
We get for every $r>0$ the following non-degeneracy property:
\begin{equation}\label{eq:ndgcy}
\sup_{Q^-_r(t_0,x_0)} w_M \ge \min \left(w_M(t_0,x_0) + \frac{1}{4n}r^2,
  A/2\right)   
\end{equation}
Combined with the Harnack-type inequality (\ref{eq:ns1}) for some radius
$r'<r$, we obtain the following bound from below:
\begin{equation}\label{eq:1.12}
\tilde{w}_M(s,x) \ge \alpha'>0 \quad \mbox{for}\quad s\le s_1 <0
\end{equation}
for some constants $\alpha'$ and $s_1$
and every large $M$.\\
%%%%%%%%%%%%%%%%%%%%%%%%%%%%%%%%%%%%%%%%%%%%%%%%%%%%%%%%%%%%%%%%%%%%%%%%%%%%%%%%
\noindent {\bf Step 6: Further estimates on the time derivative of the limit solution}\\
By the bound from above in Step 3, we obtain that
\begin{equation}\label{w_bound} |w(t,x)| \le C_1 + C_2 (|t| + |x|)\; ,\end{equation}
where $C_1$ and $C_2$ are finite positive constants.
Let now $(t_k,x_k)\in \{ w>0\}$ be a sequence such that
$$ct_k-x_k\cdot e\to -\infty\; .$$
Then by the result in Step 5, 
\begin{equation}\label{d_bound}
d_k := \pardist((t_k,x_k), \partial \{ w>0\})
\ge c_3 \sqrt{|t_k| + |x_k|^2}\end{equation}
for some constant $c_3>0$. So 
$w$ is a solution of $\partial_t w - \Delta w = -v^0$ in $Q_{d_k}(t_k,x_k)$.
Defining $$z_k(t,x) :=  \frac{w(t_k+d_k^2 t,x_k+d_k x)}{d_k^2}\; ,$$
(\ref{w_bound}) and (\ref{d_bound}) imply that
$z_k$ is a solution of $\partial_t z_k - \Delta z_k = -v^0(x_k+d_k x)$ in $Q_1(0)$
satisfying
$$\sup_{Q_1(0)} |z_k| \le C_4\; ,$$
where $C_4$ is a constant not depending on $k$.
Consequently $\partial_t w(t_k,x_k)=\partial_t z_k(0)$
is bounded, implying that $\limsup_{(t,x)\in \{ w>0\}, ct-x\cdot e\to -\infty}|\partial_t w(t,x)|<+\infty$.
Passing if necessary to a subsequence, we obtain by the periodicity
of $v^0$ a limit $z$ satisfying
$\partial_t z - \Delta z = -\mu_0$ in $Q_1(0)$.
Moreover we infer from the fact that
$\tilde w$ is periodic in the space variables
that $z$ is constant in the space variables.
Thus $\partial_t z \equiv -\mu_0$ in $Q_1(0)$.
\\
From regularity theory of caloric functions it follows that
$$\lim_{(t,x)\in \{ w>0\}, ct-x\cdot e\to -\infty} \partial_{tt} w
=0\; .$$
But then a combination of the comparison principle (\ref{eq:td})
and of (\ref{eq:1.12}) yield
$$-\mu_0\le \partial_t w \le  0  \quad \mbox{on}\quad \R\times \R^n$$
and
$$\partial_{tt} w \ge  0  \quad \mbox{on}\quad \R\times \R^n\; .$$
This ends the proof of the Theorem \ref{pro:PTWOP}.
\section{Non-existence of pulsating waves in the case of constant initial
concentration}\label{nonexpw}
We consider solutions $u$ of the one-phase limit problem with
constant initial concentration in {\em any finite dimension}, i.e.
\begin{equation}\label{eq:TW}
\partial_t u -\partial_t \chi_{\left\{u\ge 0\right\}}= \Delta u
\quad\mbox{in}\quad \R \times \R^n\; ,
\end{equation}
and prove that $u$ cannot be a non-trivial pulsating wave in the sense of
(\ref{eq:uv0}), (\ref{eq:uv01}). More precisely:
\begin{theo}\label{th:LiouvillePTW}{\bf (Non-existence of pulsating
waves for constant initial concentration)}\\
Let $u$ be a solution of (\ref{eq:uv0}), (\ref{eq:uv01}) in dimension $n\ge 1$ with $v^0=
\mbox{constant}$. Then $u(t,x)=u(t-e\cdot x /c,0)$, i.e. $u$ is a planar wave.
\end{theo}
\noindent {\it Proof of Theorem \ref{th:LiouvillePTW}}\\
We set $w(t,x)=-\int_{t}^{+\infty} u(s,x)\ ds \ge 0$.
From the 
proof of Theorem \ref{th:PTW} we know that
$w>0$ if and only if $u<0$.
As $\partial_t u \ge 0$ we obtain
$$\partial_t w = \Delta w -v^0 \chi_{\left\{w >0\right\}}\; ,$$
and $w$ satisfies (\ref{eq:w}), (\ref{eq:w1}).
For any $\xi\in\R^n$, we define the ``tangential difference''
$$z^\xi(t,x)=w(t-\frac{e\cdot \xi}{c},x-\xi)-w(t,x)$$
which satisfies
\begin{equation}\label{eq:Lc}
(\partial_t - \Delta)z^\xi = -az^\xi,
\end{equation}
$$\mbox{where}\quad 0\le a =
\left\{\begin{array}{ll}
0 & \quad \mbox{if}\quad z^\xi(t,x)=0\\
\\
\displaystyle{v^0\left(\frac{\chi_{\left\{{w}(t-\frac{e\cdot \xi}{c},x-\xi)>0\right\}}-\chi_{\left\{{w}(t,x)>0\right\}}}{{w}(t-\frac{e\cdot \xi}{c},x-\xi)-{w}(t,x)}\right)} & \quad \mbox{if}\quad z^\xi(t,x)\not= 0\; .
\end{array}\right.   
$$
From (\ref{eq:w}) ,(\ref{eq:w1}) and the definition of $z^\xi$ we infer that
$$|\partial_t z^\xi| \le 2 \mu_0 =2v^0 \textrm{ in } \R^{n+1}\; ,$$
$$\partial_t z^\xi(t,x) \longrightarrow 0 \textrm{ uniformly in }t,x,\xi \textrm{ as}\quad ct-{e\cdot x}
\longrightarrow -\infty\; .$$
Moreover (\ref{eq:w1}) and (\ref{eq:1.12}) as well as the
definition of $z^\xi$ tell us that for some $s_0\in (0,+\infty)$ not depending on $\xi$,
$$(\partial_t - \Delta) z^\xi =0 \quad \mbox{in}\quad \left\{ |ct-{e\cdot x} |> s_0\right\}\; .$$
Furthermore we obtain from the comparison principle (see Lemmata 3.2 and
3.4 in \cite{hamel}) that
\begin{equation}\label{eq:est}
|\partial_t z^\xi(t,x)|\le 2v^0 e^{ct-e\cdot x + s_0}
\end{equation}
and --- integrating this estimate for $t\in (-\infty, \frac{s_0+e\cdot x}{c})$ 
and using that $z^\xi=0$ in $t\ge \frac{s_0+e\cdot x}{c}$ ---
we obtain that $z^\xi$ is bounded on $\R^{n+1}$ by a constant not depending on $\xi$.
\\
Liouville's theorem for the heat equation implies therefore that
for each sequence $(t_m,x_m)$ such that $ct_m-{e\cdot x_m}\to -\infty$,
$z^\xi(t_m+\cdot,x_m+\cdot)$ converges locally uniformly in $\R^{n+1}$ 
(and uniformly with respect to $\xi$)
to a constant
$K$ depending on the choice of $\xi$ and the sequence $(t_m,x_m)$.
As we know that $\int_{x + [0,1)^n} z^\xi(t,y)\ dy = 0$ for every $(t,x)\in \R^{n+1}$
(see (\ref{eq:w1})),
it follows that $K=0$ and that
$$z^\xi(t+\cdot,x+\cdot)\to 0\textrm{ locally uniformly in }  \R^{n+1}
\textrm{ as } ct-{e\cdot x}\to -\infty\; ;$$
the convergence is also uniform with respect to $\xi$.\\
Finally we define
$$ \eta(t,x) := \sup_{\xi \in \R^n} |z^\xi(t,x)|\; .$$
The function $\eta$ is by (\ref{eq:Lc}) a bounded subcaloric function.
Moreover, by construction,
$$ \partial_y  \eta(t-\frac{e\cdot y}{c},x-y) \equiv 0\; .$$
But then $\eta(t,x)=f(ct-e\cdot x)$,
$cf'-f''\le 0$ in $\R$, $f\in W^{1,1}_{\rm loc}(\R)$, $\lim_{s\to -\infty} f(s)=\lim_{s\to +\infty} f(s)=0$ 
and $f$ is bounded from above, implying that
$f=0$, that $\eta\equiv 0$ and that $w(t-\frac{e\cdot \xi}{c},x-\xi)-w(t,x)=0$ for every $t\in \R$ and
$x,\xi\in \R^n$.
We obtain the corresponding result for $u$.
\section{Conclusions}\label{conclusions}
In this section we try to take a conclusion of the previous sections.
\\
Let us consider a sequence of solutions $(u_\varepsilon,v_\varepsilon)$ of the $\varepsilon$-problem
(\ref{sys}) 
satisfying for example the assumptions in the time-increasing case
of Theorem \ref{main},
and suppose that they are getting closer and closer to one-phase pulsating waves as $\varepsilon\to 0$,
i.e. for some $t_\varepsilon\to +\infty$,
\begin{equation}
\left\{\begin{array}{lcl}
\displaystyle{u_\varepsilon(t,x+k)=o(1)+u_\varepsilon(t-\frac{e\cdot k}{c},x)
\quad\mbox{for every}\quad
k\in \Z^n,\> (t,x)\in (t_\varepsilon,+\infty)\times \R^n}\\
\\
u_\varepsilon(t,x)=o(1) \quad\mbox{for}\quad x\cdot e -ct \le -t_\varepsilon\quad\mbox{and}\quad
\displaystyle{\lim_{\varepsilon\to 0} \limsup_{x\cdot e -ct +t_\varepsilon \to +\infty } u_\varepsilon(t,x) = -\mu_0}\; .
\end{array}\right.
\end{equation}
Translating $t_\varepsilon$ to $0$, we obtain by Theorem \ref{main} a sequence $(u_\varepsilon(t_\varepsilon+t,x),v_\varepsilon(t_\varepsilon+t,x))$
that is locally compact in $L^1$ and $u_\varepsilon$ converges to $u$
satisfying all assumptions in Theorem \ref{th:PTW} except the assumptions
for $v^0(x)$.
Let us assume $v^0\ge c >0$, i.e. there has been enough fuel everywhere initially,
and let us show that $v^0$ is $\Z^n$-periodic: for the integrated function $w$
as in section \ref{expw} and $k\in \Z^n$,
$$0\; = \;  (\partial_t-\Delta)(w(t,x+k)-w(t-\frac{e\cdot k}{c},x))$$
$$= v^0(x)\chi_{\{w(t-\frac{e\cdot k}{c},x)>0\}}-v^0(x+k)\chi_{\{w(t,x+k)>0\}}
= (v^0(x)-v^0(x+k))\chi_{\{w(t,x+k)>0\}}$$
and therefore, choosing $-t$ large, $v^0(x)-v^0(x+k)=0$.
\\
But then Theorem \ref{th:LiouvillePTW} tells us that $u$ must either be one of the
pulsating waves constructed in section \ref{expw}, or a planar wave.
\\
So our results when combined, strongly suggest a relation of the numerically
observed pulsating waves for dimensions $n \ge 1$ and the pulsating waves
constructed in section \ref{expw}.
Of course there remain possible alternatives, for examples our results do not
exclude the non-existence of pulsating waves for small $\varepsilon$, or the strong
convergence of pulsating waves to planar waves from far away, or the existence
of non-trivial true two-phase pulsating waves in the case of constant $v^0$.
We leave these possibilities open to mathematical discussion. However in communication with
mathematicians doing numerics for pulsating waves of the very systems
mentioned in the introduction the above alternatives have been considered
unlikely.
\\
Let us also mention that the chaos reported in \cite{chaos} for finite activation
energy presents no contradiction to our results: we never claimed that
the dynamics {\em around} the pulsating waves is simple.
\\
On the other hand, our result -- proving rigorously existence of pulsating waves
for spatially inhomogeneous coefficients and non-existence for constant
coefficients in any dimension -- suggests that numerics -- introducing {\em nolens volens}
spatially inhomogeneous coefficients -- makes the limit problem more
unstable (see also the Appendix), and that grid information etc. {\em will} strongly influence the numerically
obtained pulsating waves unless appropriate counter-measures are taken.  
\section{Open questions}
The most pressing task is of course to 
study for space dimension $n\ge 2$ the existence or non-existence of ``peaking''
of the solution in the negative phase.
A related question
is whether $(u_\varepsilon)_{\varepsilon \in (0,1)}$ is bounded in $L^\infty$
in the case of uniformly bounded initial data. Although this
seems obvious, it is not clear how to prevent concentration
close to the interface. \\
Another challenge is to use the information on the limit problem
gained in the present paper to construct pulsating waves
for the $\varepsilon$-problem.
\\
Uniqueness for the limit problem
(the irreversible Stefan problem for supercooled water) 
in general seems unlikely. One might however ask whether
time-global uniqueness holds in the case that $u$ is strictly increasing
in the $x_1$-direction. By the result in
\cite{hele} for the ill-posed Hele-Shaw problem, 
time-local uniqueness is likely to be true here, too.
\section{Appendix: Formal stability in the case of
one space dimension and constant initial concentration}\label{stability}
We consider solutions $u$ of the one-phase limit problem with
constant initial concentration in one space dimension, i.e.
\begin{equation}\label{seq:TW}
\partial_t u -\partial_t \chi_{\left\{u\ge 0\right\}}= \partial_{xx} u
\quad\mbox{in}\quad \R \times \R
\end{equation}
that are close to the traveling wave solution
\begin{equation}\label{eq:TW0}
\bar u(t,x)= - \max\left(1-e^{-c(x-ct)},0\right)   
\end{equation}
moving with velocity $c>0$.
\begin{pro}\label{pro:stab}{\bf (Formal linear stability of one-dimensional
    traveling waves)}\\
The traveling wave $\bar u$ given by (\ref{eq:TW0}) is formally linearly
stable with respect to equation (\ref{seq:TW}).
\end{pro}
\begin{rem}
In higher dimensions this result is no longer true. It is well known that
a fingering instability occurs. However the pulsation phenomenon
with which we are concerned in the present paper appears
already in dimension $1$. 
\end{rem}
\noindent {\sl Formal proof of Proposition \ref{pro:stab}:}\\
Let us consider solutions $u$ of (\ref{seq:TW}) satisfying
\begin{equation}\label{eq:u(t,x)}
\left\{\begin{array}{l}
u(t,x)=0 \quad\mbox{for}\quad x\le s(t)\; ,\\
\\
u(t,x)<0  \quad\mbox{for}\quad x> s(t)\; ,\\
\\
u(t,x)  \longrightarrow -1 \quad\mbox{as}\quad x-s(t) \to +\infty\; ,\\
\\
s'(t)=-\partial_x u(t,s(t)+0^+) \ge 0\; .
\end{array}\right.   
\end{equation}
Let us remark that a simple analysis shows that we do not have a comparison
principle for solutions of (\ref{eq:u(t,x)}).\\
In order to analyze the stability we transform (\ref{eq:u(t,x)}) by
$$v(t,y):=u(t,y+s(t))\; ;$$
$v$ satisfies $v(t,y)=0$ for $y\le 0$ and
\begin{equation}\label{eq:v(t,y)}
\left\{\begin{array}{l}
v(t,0)= 0 \; ,\\
\\
v(t,y) <0 \quad\mbox{for}\quad y> 0\; ,\\
\\
v(t,y)  \longrightarrow -1 \quad\mbox{as}\quad y \to +\infty\; ,\\
\\
\partial_t v = \partial_{yy}v +s'(t) \partial_y v \quad\mbox{on}\quad \R
\times (0,+\infty )\; ,\\
\\
s'(t)=-\partial_y v(t,0^+)\; .
\end{array}\right.   
\end{equation}
We now consider for $t>0$ a perturbation of the traveling wave $$\bar v(y)=-
\max\left(1-e^{-cy},0\right)$$ with velocity $c>0$.
In the formal expansion
$$\left\{\begin{array}{l}
s(t)=ct+ \varepsilon \gamma(t) + 0(\varepsilon^2)\; ,\\
\\
v=\bar v +\varepsilon w + 0(\varepsilon^2)\; ,
\end{array}\right.$$
the first order terms $w(t,y),\gamma(t)$ formally satisfy
$$\left\{\begin{array}{l}
w(t,0)= 0 \; ,\\
\\
w(t,y)  \longrightarrow 0 \quad\mbox{as}\quad y \to +\infty\; ,\\
\\
\partial_t w = \partial_{yy}w + c\partial_y w  + \gamma'(t) \partial_y \bar v\quad\mbox{in}\quad (0,+\infty )
\times (0,+\infty )\; ,\\
\\
\gamma'(t)=-\partial_y w(t,0^+)\; .
\end{array}\right.$$
Let us look for solutions of the form
$$\left\{\begin{array}{lcl}
w(t,y)=e^{\lambda t}W(y)   \; ,\\
\\
\gamma'(t)=e^{\lambda t}\; ,
\end{array}\right.$$
where $\mbox{Re }(\lambda)\ge 0$.
We obtain
$$W''+c W'-\lambda W = ce^{-cy}\; ,$$
i.e.
$$W(y)=-\frac{c}{\lambda}e^{-cy} +\sum_{\pm} A_{\pm} e^{\mu_\pm y}\; ,$$
where
$$\mu_\pm = -c/2 \pm \sqrt{c^2/4 +\lambda} \quad\mbox{and}\quad \mbox{Re
}\left(\sqrt{c^2/4 +\lambda}\right)>c/2\; .$$
The function $W$ can only be bounded if
$A_+=0$ and, by $W(0)=0$,
$$W(y)=-\frac{c}{\lambda}\left(e^{-cy} -e^{\mu_- y}\right)\; .$$
Finally the relation $\gamma'(t)=-\partial_y w(t,0)$ implies
$$1= \frac{c}{\lambda} \left(-c -\mu_-\right)\; .$$
The unique solution of this equation is $\lambda=0$. 
Thus we formally proved stability of traveling waves.
%%%%%%%%%%%%%%%%%%%%%%%%%%%%%%%%%%%%%%%%%%%%%%%%%%%%%%%%%%%%%%%%%%%%%%%%%%%%%
\newline
% \begin{pro}\label{pro:noex}{\bf (Formal non-existence of one-dimensional
% pulsating waves)}\\
% Formally, solutions $u$ of equation (\ref{eq:TW}) that are not
% traveling waves, but satisfy for some period $T>0$ and velocity $c>0$
% \begin{equation}\label{eq:PW}
% u(t+T,x)=u(t,x-cT)\; ,   
% \end{equation}
% do not exist.
% \end{pro}
% \noindent {\sl Formal proof of Proposition \ref{pro:noex}:}\\
% Let us once more use the function $v(t,y)=u(t,y+s(t))$ satisfying
% (\ref{eq:v(t,y)}).
% Now assume that $v(t,\cdot)$ is the following unilateral Laplace transform of some 
% measure $\mu(t,\cdot)$:
% $$v(t,y)=\int_{-\infty }^0 \mu(t,\beta) e^{\beta y}\ d\beta$$
% We see that $\mu$ formally satisfies
% \begin{equation}\label{eq:f}
% \partial_t \mu=\beta^2 \mu + s'(t) \beta \mu   
% \end{equation}
% which yields by integration
% \begin{equation}\label{eq:beta}
% \mu(t,\beta)=\mu(0,\beta)e^{\beta\left(s(t)+\beta t -s(0)\right)}\; .   
% \end{equation}
% Assuming (\ref{eq:PW}), we deduce from (\ref{eq:u(t,x)}) that
% $$s(t+T)=s(t)+cT\; .$$
% Moreover we infer from (\ref{eq:v(t,y)}) that
% $$v(t+T,y)=v(t,y)\; .$$
% Formally it follows that
% $$\mu(t+T,\beta)=\mu(t,\beta)\; ,$$
% whence 
% we conclude by (\ref{eq:beta}) that either $\beta =0$ or
% $\beta=-c$. Thus the measure $\mu(t,\cdot)$ must be concentrated
% at those two points.
% We get
% $$v(t,y)=C_1 +C_2 e^{-cy}\; .$$
% By (\ref{eq:v(t,y)}), $C_1=-1$ and $C_2=1$, implying -- at least formally -- that
% the only
% pulsating wave traveling at velocity $c>0$ is the traveling wave.
{\bf Acknowledgment:}
We thank Steffen Heinze, Danielle Hilhorst, Stephan Luckhaus, Mayan Mimura, Stefan M\"uller 
and Juan J.L. Vel\'azquez for discussions.
  
%%%%%%%%%%%%%%%%%%%%%%%%%%%%%%%%%%%%%%%%%%%%%%%%%%%%%%%%%%%%%%%%%%%%%%%%%%%%%%%%%

\bibliographystyle{plain}
\bibliography{monneauweiss_shs061130.bib}
\end{document}